\documentclass[11pt]{article}
\usepackage{newlfont,amsfonts,amssymb, oldgerm,amsmath,amsthm,amsgen,amscd,
}
\usepackage{epic}
\usepackage{eepic}
\usepackage{rotating}
\usepackage[dvips]{epsfig}
\DeclareMathOperator{\Sig}{Sig}
\DeclareMathOperator{\codim}{codim}

\begin{document}
\theoremstyle{definition}
\newtheorem*{ex}{Example}
\newtheorem*{theo*}{Theorem}
\newtheorem{de}{Definition}
\newtheorem{bem}{Remark}
\newtheorem{bez}{Notation}
\theoremstyle{plain}
\newtheorem*{con}{Conjecture}
\newtheorem{lem}{Lemma}
\newtheorem{satz}{Proposition}
\newtheorem{folg}{Corollary}
\newtheorem{theo}{Theorem}
\newtheorem{fact}{Fact}

\newcommand{\lan}{\mathfrak{n}}
\newcommand{\lah}{\mathfrak{h}}
\newcommand{\laz}{\mathfrak{z}}

\newcommand{\rrn}{\mathbb{R}^n}
\newcommand{\laso}{\mathfrak{so}}
\newcommand{\lason}{\mathfrak{so}(n)}
\newcommand{\lagl}{\mathfrak{gl}}
\newcommand{\lasl}{\mathfrak{sl}}
\newcommand{\lasp}{\mathfrak{sp}}
\newcommand{\lasu}{\mathfrak{su}}

\newcommand{\rr}{\mathbb{R}}
\newcommand{\ccc}{\mathbb{C}}
\newcommand{\bs}{\begin{satz}}
\newcommand{\es}{\end{satz}}
\newcommand{\btheo}{\begin{theo}}
\newcommand{\etheo}{\end{theo}}
\newcommand{\bfolg}{\begin{folg}}
\newcommand{\efolg}{\end{folg}}
\newcommand{\blem}{\begin{lem}}
\newcommand{\elem}{\end{lem}}
\newcommand{\bprf}{\begin{proof}}
\newcommand{\eprf}{\end{proof}}
\newcommand{\bd}{\begin{displaymath}}
\newcommand{\ed}{\end{displaymath}}
\newcommand{\be}{\begin{eqnarray*}}
\newcommand{\ee}{\end{eqnarray*}}
\newcommand{\eeqa}{\end{eqnarray}}
\newcommand{\beqa}{\begin{eqnarray}}
\newcommand{\bi}{\begin{itemize}}
\newcommand{\ei}{\end{itemize}}
\newcommand{\bnum}{\begin{enumerate}}
\newcommand{\enum}{\end{enumerate}}
\newcommand{\la}{\langle}
\newcommand{\ra}{\rangle}
\newcommand{\ve}{\varepsilon}
\newcommand{\vp}{\varphi}
\newcommand{\lra}{\longrightarrow}
\newcommand{\Lra}{\Leftrightarrow}
\newcommand{\Ra}{\Rightarrow}
\newcommand{\sub}{\subset}
\newcommand{\ems}{\emptyset}
\newcommand{\sms}{\setminus}
\newcommand{\ints}{\int\limits}
\newcommand{\sums}{\sum\limits}
\newcommand{\lims}{\lim\limits}
\newcommand{\bcup}{\bigcup\limits}
\newcommand{\bcap}{\bigcap\limits}
\newcommand{\beq}{\begin{equation}}
\newcommand{\eeq}{\end{equation}}
\newcommand{\mf}{\mathfrak}
\newcommand{\lag}{\mathfrak{g}}
\newcommand{\+}{\oplus}
\newcommand{\einhalb}{\frac{1}{2}}

\newcommand{\C}{\ensuremath{\mathbb{C}}}
\newcommand{\R}{\ensuremath{\mathbb{R}}}
\newcommand{\Z}{\ensuremath{\mathbb{Z}}}
\newcommand{\N}{\ensuremath{\mathbb{N}}}
\newcommand{\K}{\ensuremath{\mathbb{K}}}
\newcommand{\HH}{\ensuremath{\mathbb{H}}}
\newcommand{\II}{\ensuremath{\mathbb{I}}}
\newcommand{\JJ}{\ensuremath{\mathbb{J}}}
%
\newcommand{\G}{\ensuremath{\mathcal{G}}}
\newcommand{\I}{\ensuremath{\mathcal{I}}}
\newcommand{\RR}{\ensuremath{\mathcal{R}}}
\newcommand{\T}{\ensuremath{\mathcal{T}}}
\newcommand{\X}{\ensuremath{\mathcal{X}}}
\newcommand{\Kill}{\mathcal{B}^{Killing}}
\newcommand{\B}{\ensuremath{\mathcal{B}}}
\newcommand{\Zent}{\ensuremath{\mathcal{Z}}}
%

\newcommand{\g}{\ensuremath{\mathfrak{g}}}
\renewcommand{\k}{\ensuremath{\mathfrak{k}}}
\newcommand{\h}{\ensuremath{\mathfrak{h}}}
\renewcommand{\i}{\ensuremath{\mathfrak{i}}}
\newcommand{\m}{\ensuremath{\mathfrak{m}}}
\newcommand{\n}{\ensuremath{\mathfrak{n}}}
\newcommand{\p}{\ensuremath{\mathfrak{p}}}
\newcommand{\s}{\ensuremath{\mathfrak{s}}}
\newcommand{\z}{\ensuremath{\mathfrak{z}}}
\renewcommand{\r}{\ensuremath{\mathfrak{r}}}
\newcommand{\ad}{\ensuremath{\mathfrak{ad}}}
\newcommand{\gl}{\ensuremath{\mathfrak{gl}}}
\newcommand{\so}{\ensuremath{\mathfrak{so}}}
\renewcommand{\o}{\ensuremath{\mathfrak{o}}}
\newcommand{\su}{\ensuremath{\mathfrak{su}}}
\renewcommand{\sl}{\ensuremath{\mathfrak{sl}}}
\newcommand{\spin}{\ensuremath{\mathfrak{spin}}}
\newcommand{\rad}{\ensuremath{\mathfrak{rad}}}
\newcommand{\iso}{\ensuremath{\mathfrak{is\mspace{-2mu} o}}}
\newcommand{\hol}{\ensuremath{\mathfrak{hol}}}
\newcommand{\aut}{\ensuremath{\mathfrak{aut}}}
\newcommand{\mb}[1]{{\mathversion{bold}{\bf #1}}}
\newcommand{\ol}{\overline}
\newcommand{\rca}[2]{ \renewcommand{\arraystretch}{#1}{#2} }
\newcommand{\wt}{\widetilde}
\newcommand{\bdm}{\begin{displaymath}}
\newcommand{\edm}{\end{displaymath}}
\newcommand{\btab}{\begin{tabular}}
\newcommand{\etab}{\end{tabular}}
\newcommand{\ba}{\begin{array}}
\newcommand{\ea}{\end{array}}
\newcommand{\bsm}{\begin{smallmatrix}}
\newcommand{\esm}{\end{smallmatrix}}
\newcommand{\bal}{\begin{align*}}
\newcommand{\eal}{\end{align*}}
\newcommand{\comment}[1]{{\it[Comment: #1 ].}}
\newcommand{\menge}[2]{ \big\{\, #1 \,\big|\, #2 \big\} }
\newcommand{\ds}{\displaystyle}
\newcommand{\sst}{\scriptstyle}
\newcommand{\lp}{\left(}
\newcommand{\rp}{\right)}

\bibliographystyle{alpha}

\title{Irreducibly acting subgroups of $Gl(n,\rr)$}

\author{Antonio J. Di Scala\thanks{ Dipartimento di Matematica, Politecnico di Torino, 
Corso Duca degli Abruzzi, 24, 
10129 Torino, Italy. Email: {\tt antonio.discala@polito.it}},
Thomas Leistner\thanks{School of Mathematical Sciences, University of Adelaide,  SA 5005, Australia.
 Email: {\tt tleistne@maths.adelaide.edu.au}
}, and Thomas Neukirchner\thanks{{ Institut f\"{u}r Mathematik, Humboldt-Universit\"{a}t Berlin, Unter den Linden 6, D-10099 Berlin, Germany.}
 Email: {\tt neukirch@mathematik.hu-berlin.de}}}

\date{ }
\maketitle
\begin{abstract}
In this note we prove the following three algebraic facts which have applications in the theory of holonomy groups and homogeneous spaces: Any irreducibly acting connected subgroup $G \subset Gl(n,\rr)$ is closed. Moreover, if $G$ admits an invariant bilinear form of Lorentzian signature, $G$ is maximal, i.e.
it is conjugated to $SO(1,n-1)_0$. Finally we calculate the vector space of $G$-invariant symmetric bilinear forms, show that it is  at most $3$-dimensional, and determine
the maximal stabilizers for each dimension.
\\[.1cm]
\end{abstract}


\section{Background,  results and applications}  
This article deals with three algebraic questions which are motivated by problems in holonomy theory of affine or semi-Riemannian manifolds and in the theory of homogeneous spaces. The three problems are the following: Are holonomy groups closed?  What are special holonomy groups of Lorentzian manifolds? And finally, how many $G$-invariant bilinear forms exist on a homogeneous space $G/H$?

Regarding the first question, the first thing to observe is that the answer is `no' in general,  because due to a result of   \cite{hano-ozeki} any linear Lie group can be realised as a holonomy group of a linear connection. Secondly one finds that the answer is `yes' if one restricts it to the holonomy group of the Levi-Civita connection of a Riemannian manifold, because 
 any {\em irreducibly acting}, connected subgroup of $O(n)$ is closed
(see for example \cite[Appendix 5]{KobayashiN1}), and thus,  by
 the de Rham decomposition  theorem  \cite{derham52} the
holonomy group is a direct product of closed ones.
But this result can not be generalised to the Levi-Civita connection of semi-Riemannian manifolds (see \cite{wu67} for pseudo-Riemannian manifolds and
\cite{bb-ike93} for Lorentzian manifolds with non-closed holonomy groups) and not even to torsion-free connections (see examples in \cite{hano-ozeki}). In all of these examples the holonomy group does not act irreducibly.
So it arises the question if all irreducibly acting subgroups of $Gl(n,\rr)$ are closed? There is a note in \cite[Note 10.50, page 290]{besse87}  that the following theorem is true, proved in \cite{wakakuwa71}. 
\btheo\label{p1}
Any irreducibly acting connected Lie subgroup of $Gl(n,\rr)$ is closed in $Gl(n,\rr)$.
\etheo
In this note firstly we 
will prove this theorem independently on the proof in \cite{wakakuwa71}.
Regarding holonomy groups of linear connections this has the following consequence.
\bfolg\label{folge1}
If the restricted holonomy group of a linear connection acts irreducibly, then it is closed. Furthermore,  the restricted holonomy group of a semi-Riemannian manifold is closed if it acts completely reducibly.
\efolg
The second statement follows from the de Rham/Wu decomposition theorem \cite{wu64} and another theorem in \cite{wu67} (see below, last paragraph in section \ref{sec2}).
Our proof of  Theorem \ref{p1} uses a theorem of Yosida \cite{yosida37} and Malcev \cite{malcev45}  and a very explicit description of the center of $G$. This description gives us two corollaries, the first of which will be useful in the proof of Theorem \ref{p2}.
\bfolg\label{folge3}
Let $G\subset O(p,q)$ be a connected Lie subgroup of $O(p,q)$ which acts irreducibly. If $G$ is not semisimple, then $p$ and $q$ are even and $G$ is a subgroup of $ U(p/2, q/2)$ with center $U(1)$. 
\efolg  

Applying this to the spin representations of orthogonal groups we get:
\bfolg\label{folge2}
Let $G\subset SO_0(p,q)$ be a connected Lie subgroup which acts irreducibly and $\tilde{G}\subset Spin(p,q)$ its lift into the spin group. If  the spin representation of $\tilde{G}$ admits a trivial subrepresentation, then $G$ is semisimple. 
\efolg
Also this corollary    has applications to geometric problems. The first application is a well known fact. If the holonomy group of a semi-Riemannian manifold acts irreducibly and  has a center, the manifold cannot admit parallel spinors. It was obtained by the classification of irreducible holonomy groups of semi-Riemannian manifolds with parallel spinors (\cite{wang89} for Riemannian manifolds, and \cite{baumkath99} for pseudo-Riemannian manifolds). 
But furthermore it gives results in any case where the holonomy group has a irreducibly acting component on which the existence of parallel spinors depends, as it is the case for indecomposable, non-irreducible Lorentzian manifolds, see \cite{leistner01}. 

\bigskip

Regarding the second problem which special Lorentzian holonomy groups might exist, one distinguishes between the irreducible and the indecomposable, non-irreducible case. While in the latter case there are several possibilities (for a classification see \cite{bb-ike93}, \cite{leistner02}, \cite{leistner03}, \cite{leistner03b}, 
and \cite{galaev05}), for the irreducible case the situation is very limited.
The irreducible holonomy groups of semi-Riemannian manifolds were determined by M. Berger in \cite{berger55} and \cite{berger57}. In Riemannian and many other signatures  the list
depends essentially on the property of being a \emph{holonomy group}, whereas in the Lorentzian case it turns out that \emph{irreducibility}
is sufficient to determine the group.
\btheo\label{p2}
$SO(1,n)_0$ is the only connected Lie subgroup of $O(1,n)$ which acts irreducibly.
\etheo
The consequence for irreducible Lorentzian holonomy groups follows immediately.
\bfolg
If the restricted holonomy group of a Lorentzian manifold acts irreducibly, then it is equal to 
$SO(1,n)_0$.
\efolg
A direct and geometric proof of Theorem \ref{p2}  was given in \cite{olmos-discala01}. In Section \ref{sec3}  we will give a short proof of Theorem \ref{p2} based on a theorem of Karpelevich \cite{Ka} and Mostow \cite{Mo}.

\bigskip

The result of the last section is motivated by the geometric problem of 
describing the space of metrics or symplectic forms on  a homogeneous space $G/H$ which are invariant under $G$. Any $G$-invariant
metric or symplectic form corresponds to a non-degenerate bilinear form on $\g/\h$  which is invariant under the
linear isotropy representation $Ad_G(H) \subset Gl(\g/\h)$. In our context $Ad_G(H)$ is assumed to act irreducibly. 
This is a special case of the following algebraic problem: Given an irreducibly acting Lie subgroup $G\subset Gl(n,\rr)$, what is the dimension of the space of $G$-invariant bilinear forms on $\rrn$. We prove the following statement.
\btheo\label{p3}
Let $G$ be an irreducibly acting subgroup of $Gl(n,\rr)$. The space of $G$-invariant symmetric bilinear forms which are not of neutral signature $(p,p)$ is at most one-dimensional.
Moreover, the space of invariant symmetric bilinear forms is at most three-dimensional.
\etheo
We will describe all possible cases for the dimension of the space of (skew-) symmetric bilinear forms
and determine the maximal subgroup which fixes these bilinear forms.
 
\bigskip
 
We should point out that many results in this paper rely on the classification of $G$-invariant endomorphisms for $G\subset Gl(n,\rr)$. This classification follows from Schur's lemma and the classification of associative division algebras by Frobenius, but we will give an elementary proof of it in Section \ref{sec1}.

\section{The algebra of invariant endomorphisms}
\label{sec1}

The results of this paper are mainly based on a description of the algebra of endomorphisms which are invariant under an irreducibly acting subgroup $G\subset Gl(n,\rr)$. If $G$ is a Lie group and $V$ and $W$ two (real or complex) $G$-modules  the algebra of invariant homomorphism is defined as
\[Hom_G(V,W)\ :=\ \{X\in Hom(V,W)\ |\ A\circ X =X\circ A\ \text{for all } A\in G\}.\]
Now, Schur's lemma says that $Hom_G(V,W)\subset Iso(V,W)\cup \{0\}$, and furthermore, if $V=W$ is complex, then $End_G(V)=\ccc\cdot Id$. In any case, it implies that 
$Hom_G(V,W)$ is a real associative division algebra, and thus by their classification of Frobenius (1878) it  is isomorphic to the algebra of real numbers $\R$,
complex numbers $\C$ or quaternions $\HH$ (see e.g. \cite{Palais68}). We are interested in the description of $End_G(V)$ where $V$ is a real vector space, and in this section we will recall some facts about real irreducible representations  which provide an elementary proof of this result.

\bigskip

Suppose that $G$ is a real Lie group and $V$ a real irreducible module. Then there are two cases which can occur for the complexified $G$-module $V^\ccc$. The first case is that $W:=V^\ccc$ is still irreducible. In this case $V$ or $W$ is called  {\em of real type}. One should remark that, if $W$ is a complex irreducible $G$-module then its reellification $W_\rr$ is a reducible $G$-module with invariant real subspace $V$ if and only if $W$ is the complexification of the real irreducible $G$-module $V$. 

In the other and more complicated case, regarding the application of Schur's lemma,  $V^\ccc$ is a reducible $G$-module. In this case $V^\ccc$ splits into two irreducible $G$-modules, 
\[V^\ccc=W\+\ol{W}.\]
In fact, if $W$ is an invariant complex subspace of $V^\ccc$ then $\overline{W}$, defined by the conjugation with respect to $V\subset V^\ccc$ is  invariant too and the conjugate module. Furthermore, the spaces $W+\overline{W}$ and $W\cap\overline{W}$
are invariant and equal to their conjugation.
Hence they are complexifications of  real vector spaces,
i.e. $ W+\overline{W}=V_1^\ccc$ and $W\cap\overline{W}=V_2^\ccc$. Of course $V_1$ and $V_2$ are invariant subspaces of $V$
and thus $V_1=V$ and $V_2= \{ 0 \} $. The same argument ensures the irreducibility of $W$.

Now, since $W\cap\ol{W}=\{0\}$, the mapping
 $\psi: W_\rr\ni v\mapsto \einhalb (v+\overline{v})\in V$ is an isomorphism of real vector spaces yielding the identification
 \beq \label{psi}
 W_\rr\ \stackrel{\psi}{\simeq}\  V\ \stackrel{\psi}{\simeq}\  \ol{W}_\rr\eeq
 of real $G$-modules. 
   In this case $V$, respectively $W$, are called {\em of complex type}, and again we have that a complex module $W$ has an irreducible reellification $V=W_\rr$ if and only if $V^\ccc=W\+\ol{W}$ is reducible.

%


Now we are able to describe the algebra of invariant endomorphisms of a real irreducible $G$-module $V$.
\begin{satz}\label{End_G(V)}
Let $G$ a Lie group and $V$ a real irreducible $G$-module.
Then $End_G(V)$ is isomorphic to one of the real algebras $\R$, $\C$ or $\HH$.
\end{satz}

\begin{proof} As above we  consider two cases. Firstly assume that $V^\ccc$ is irreducible which ensures that $End_G(V^\ccc)=\ccc\cdot Id$ by Schur's lemma.
Hence, if $A\in End_G(V)$, its complexification $A^\ccc\in End_G(V^\ccc)$ is given by $A^\ccc=\lambda\cdot Id$ with $\lambda\in \ccc$. Since $A^\ccc$ leaves $V$ invariant and $V$ is invariant under conjugation we get for $v=\ol{v}\in V$ that
\[\ol{\lambda}\ v\ =\ \ol{\lambda}\ \ol{v}\ =\ \ol{A^\ccc v}\ =\ A^\ccc v\ =\ \lambda v,\]
i.e. $\lambda\in \rr$. $A=A^\ccc |_V$ gives that $End_G(V)=\rr\cdot Id$.

For the second case we have to assume that $V^\ccc$ is reducible, i.e. by the above
$V^\ccc=W\+\ol{W}$, $V\simeq W_\rr$ and thus $End_G(V)=End_G(W_\rr)$.
Now any real endomorphism on $W_{\R}$ decomposes uniquely into a complex linear and complex anti-linear part:
\[
\ba{rcccc} End(W_{\R}) &\simeq&  End(W) &\oplus& Hom(W,\ol{W})\\
A &=& \frac{1}{2}(A+iAi) &+&  \frac{1}{2}(A-iAi). \ea
\]
This decomposition descends to $End_{G}(W_{\R})$:
\[
End_{G}(W_{\R}) \simeq  End_{G}(W) \oplus {Hom_{G}(W,\ol{W})}.
\]
Now
Schur's lemma implies that $End_{G}(W)= \C \cdot Id$ and, since both, $W$ and $\ol{W}$ are irreducible, that $Hom_G(W,\ol{W})\subset Iso(W,\ol{W})\cup \{0\}$.

If $Hom_G(W,\ol{W})=\{0\}$ we get immediately 
\[\ccc\cdot Id\ =\ End_G(W_\rr)\ \simeq\  End_G(V)\ =\  span_{\R}\{Id,I\},\]
with the complex structure $I:=\psi \circ (i\cdot Id) \circ \psi^{-1}$ where $\psi$ is defined in (\ref{psi}).

Otherwise consider a non-zero $j \in Hom_G(W,\ol{W})$ which is an isomorphism by Schur's lemma. Then $j^2 \in End_{G}(W)$, hence $j^2=\lambda \cdot Id$ with $0\not=\lambda \in \C$. In fact, $\lambda \in \R$
since 
\[\ol{\lambda} \ j(w) \ = \ j (\lambda \ w) \ =\  j\left( j^2(w)\right) \ =\  j^2 \left(j(w)\right) \ =\  \lambda\ j(w)\]
for all $w\in W$. Finally $\lambda<0$ since otherwise $W_\rr$ would decompose
into the $G$-invariant $\pm \sqrt{\lambda}$-eigenspaces of $j_{\R}$. Thus, we may assume $j^2=-1$. For another $A\in Hom_G(W,\ol{W})$ we get
$j \circ A \in End_{G}(W)$ and therefore $j \circ A = c \cdot Id$ for some $c \in \C$.
On the other hand $j \circ (-\ol{c}j)=c \cdot Id$ and thus $A=-\ol{c}j$. Hence we obtain 
\[End_G(W_\rr)\ \simeq\  End_G(W)\+Hom_G(W,\ol{W})\ =\ \ccc\cdot Id \+\ccc\cdot j,\] 
which gives finally
\[End_{G}(V)=span_{\R}\{Id,I,J,I \circ J\} \simeq \HH,\]
with $I:=\psi \circ i \circ \psi^{-1}$ and $J:=\psi \circ j \circ \psi^{-1}$
anti-commuting complex structures.
\end{proof}

Corresponding to the structure of $End_G(V)$ the real irreducible $G$-module $V$ is said to be of \emph{real}, \emph{complex}
or \emph{quaternionic} type. This corresponds to the convention to call a complex irreducible $G$-module $W$ of {\em real type} if it is self-conjugated with respect to an anti-linear bijection $J$ with $J^2=Id$, of {\em quaternionic type} if it is self-conjugated with $J^2=-Id$ and of {\em complex type} if it is not self-conjugated.
Here is a useful consequence of the preceeding proposition.

\begin{folg}\label{id+J}
For a real irreducible $G$-module $V$ any $A \in End_G(V)$ is of the form
$
A= \alpha \,Id + \beta J$ with $\alpha,\beta \in \R
$ and 
 $J$  a $G$-invariant complex structure (depending on $A$).
\end{folg}
\begin{proof}
Although this follows directly from Proposition \ref{End_G(V)} we will give another proof which will be useful later on. Applying the Schur-lemma
we see that the minimal polynomial $\mu_A(x)$ of $A$  is irreducible over $\R$ (cf. \cite[Appendix 5, Lemma 1]{KobayashiN1}).
If $\mu_A(x)=x-\alpha$ is of degree one $0=\mu_A(A)=A-\alpha \cdot Id$.
Otherwise $\mu_A(x)=(x - \alpha )^2 + \beta^2$ is a polynomial of degree 2 with strictly positive
quadratic supplement, since $\mu_A$ is irreducible.
Thus  $J:=(A - \alpha \cdot Id)/\beta$ defines a complex structure on $V$.
\end{proof}
Finally in this section we describe the maximal representations of different types, i.e. any other irreducible representation occurs as
a sub-representation of them.

\begin{satz}\label{Class1}
Let $G \subset Gl(n,\R)$ be an irreducibly acting subgroup.
Then up to conjugation $G$ is contained in one of the following subgroups $L \subset Gl(n,\R)$:
\[
\rca{2}{\ba{|c|c|}
\hline
End_{G}(\R^n) & L \subset Gl(n,\R) \\ \hline
\R & Gl(n,\R) \\
\C & Gl(n/2,\C) \\
\HH & Gl(n/4,\HH)
\\ \hline \ea
}
\]
\end{satz}

\begin{proof}
We set $V:=\R^n$ and $\K:=End_G(V)$. Thus $V$ becomes a \emph{left} $\K$-vector space
in a natural way. In order to make it a \emph{right} $\K$-vector space we choose an
anti-automorphism $\lambda \mapsto \ol{\lambda}$ of $\K$
(i.e. $\ol{\lambda + \mu}=\ol{\lambda}+\ol{\mu}$ and $\ol{\lambda \cdot \mu}=\ol{\mu} \cdot \ol{\lambda}$).
Then $v \cdot \lambda:=\ol{\lambda}(v)$ defines a right-multiplication on $V$ with respect to the scalar field $\K$
(This is essential only in case of the non-commutative field $\K=\HH$).
The group $Gl(V,\K)$ of $\K$-linear invertible maps from $V$ into itself is by definition
the centralizer of the homothety group $H_{\K}:=\menge{\{v \mapsto v\lambda\}}{\lambda \in \K^*}$.
By choosing a $\K$-basis $\{b_i\}_{i=1}^{n/d}$, where $d=\dim_{\R}\K$, we get
$\K^{n/d}\simeq V$. Under this identification $Gl(V,\K)$ corresponds to the group
$Gl(n/d,\K)$ of invertible $(n/d \times n/d)$-matrices acting on $\K^{n/d}$ from the left.
By definition $H_{\K}$ is the centralizer of $G$ and thus $G$ is contained in $L:=Gl(V,\K)$.
As explained this yields an inclusion $G \subset Gl(n/d,\K)$. Conversely it is known that
the centralizer of $Gl(n/d,\K)$ equals $H_{\K}$, hence $End_L(\R^n)=H_{\K}$.
Finally, the embedding
$Gl(n/d,\K) \subset Gl(n,\R)$ is obtained by associating to the $\K$-basis $\{b_i\}$ the real basis
$\{b_i \lambda_k\}_{\rca{.5}{\ba{l} \sst i=1,\ldots,n/d \\ \sst k=1,\ldots,d \ea}}$,
where $\{\lambda_k\}_{k=1,\ldots,d}$ is a basis of $\K$.
\end{proof}
\begin{bem}
In the proof of this proposition we see that if the action of a group $G\subset Gl(n,\rr)$ is defined 
 by scalar multiplication from the right, the invariant endomorphism have to act from the left.  Of course, this becomes only relevant in case of $End_G (\rrn)=\HH$, and we can see this in the example of $G:=Gl(1,\HH)$: It is
\[
Gl(1,\HH)\ =\ \{R_q:\HH\rightarrow \HH \mid q\in \HH^* \text{ and }R_q(p):=p\cdot q\}\ =\ \HH^*,
\]
whereas
\be
End_{Gl(1,\HH)}(\rr^{4})&=&\{A\in Gl(4,\rr) \mid 
A(R_q(p))=R_q(A(p)\}
\\
&=&\{L_q\in Gl(4,\rr)\mid L_q(p):=q\cdot p\}
\\
&=&\HH
\ee
since $L_q\circ R_p=R_p\circ L_q$ but $R_p\circ R_q\not=R_q\circ R_p$. This gives the seemingly paradoxical situation where both, the centraliser  $Z_{Gl(4,\rr)}(G)$ and the group $G$ itself are equal to $\HH^*$, but its center $Z(G)$ which is the intersection of $G$ with its centraliser is commutative and thus equal to $\ccc^*$.
\end{bem}
\section{Irreducibly acting, connected subgroups of \mb{$ Gl(n,\rr)$}}\label{sec2}

In this section we shall give a proof of Theorem \ref{p1} by using the results of the first section and  two general results from Lie theory. First we describe the identity component of the center of an irreducibly acting Lie subgroup of $Gl(n,\rr)$. We should remark that we mean `Lie subgroup' always in the weaker sense of being a immersed submanifold but not necessarily an embedding in order to make the statement of Theorem \ref{p1} non-trivial.
\bs
\label{center}
Let $G\subset Gl(n,\rr)$ be an irreducibly acting, connected Lie subgroup, $Z(G)$ its center and $Z(G)_0$ the identity component of the center. If  $Z(G)_0$ is non-trivial, then  $Z(G)_0$ is either 
\bnum
\item[(a)]
equal to $\rr^+ Id$, or 
\item[(b)] isomorphic  to $\ccc^*=\rr^+\times S^1$, or 
\item[(c)] isomorphic 
to a one-parameter subgroup of $\ccc^*$.
\enum
Cases $(b)$ and $(c)$ can only occur if $n$ is even.
\es

 \bprf
Let $\lag\subset \lagl(n,\rr)$ be the Lie algebra of $G$ and suppose that center $\laz$ of $\lag$ is non trivial. Considering the three cases of Proposition \ref{End_G(V)} we first assume  that the representation is of real type, i.e. that $End_G(\rrn)= \rr Id$. Since $\laz\subset End_G(\rrn)$ we obtain in this case that $\laz=\rr Id$ and therefore $Z(G)_0=\exp (\laz)=\rr^+Id$.

Now suppose that $\rrn$ is a $G$-module of non-real type, i.e. $End_G(\rr^{2n})$ isomorphic to $\ccc$ or $\HH$. Again $\laz$ is an Abelian subalgebra of $End_G(\rr^{2n})$. In case $End_G(\rr^{2n})\simeq \HH \simeq \mf{u}(2)$ any maximal Abelian subalgebra is isomorphic to $\ccc$. Hence $\laz$ is isomorphic to a subalgebra of $\ccc=span_\rr (Id, J)$ where $J $ is a complex structure on $\rr^{2n}$. But $\exp tJ=(\cos t) Id + (\sin t) J$, i.e. $\exp (\rr J)\simeq S^1$. But this implies that $Z(G)_0$ is either isomorphic to $\C^*$, i.e.
\[Z(G)_0\ =\ \rr^+ Id \times \{(\cos t) Id + (\sin t) J\mid t\in \rr\}\ \simeq \ \rr^+  \times S^1\ =\  \C^*,\]
or to a one-parameter subgroup of it, i.e. 
\be
Z(G)_0&=&\exp \left( \rr \cdot \left( a Id +b J\right)\right)\\
&=& \left\{ \left(\text{e}^{at}\cdot Id\right)\circ \left( (\cos bt) Id +(\sin bt) J\right)\mid t\in \rr\right\},
\ee
 for some real constants $a$ and $b$. Of course if $a$ or $b$ are zero this is either $\rr^+$ or $S^1$, if not this is a logarithmic spiral in $\C^*$.
 \eprf

Proposition \ref{center} will be the main ingredient in our proof of Theorem \ref{p1} but it implies also Corollaries \ref{folge3} and \ref{folge2} given in the introduction.
But before we can prove these we have to recall that for
 a completely reducibly acting Lie subgroup $G\subset Gl(n,\rr)$ the center decides whether the Lie algebra is semisimple or not. This is due to a standard fact from  the theory of Lie algebras, saying that a
Lie algebra $\lag$ which admits a  completely reducible representation is reductive. Hence $\lag$ admits a Lie algebra decomposition into its  center and its derived Lie algebra,
\beq\label{reduct} \lag \ = \ \laz \+ \left[\lag,\lag\right],\eeq
the derived Lie algebra being semisimple.
A proof of this fact can be found in \cite{chevalley47}, see also \cite{bourbaki-lie1-3}. 
This means that the irreducibly acting, connected Lie subgroup in question is semisimple if the identity component of its center is trivial. 
\begin{bem}
In this context we should remark that the center of a semisimple subgroup $G\subset Gl(n,\rr)$ is finite (see e.g. \cite{goto48}):   if $G$ is semisimple,
due to Weyl's theorem it acts completely reducibly, and furthermore its elements are of determinant $1$, hence by  Schur's lemma the center of $G$ corresponds to the $n_k$-th roots of $1$ where $n_k$ are the dimensions of the irreducible subspaces.
\end{bem}

For verifying Corollary \ref{folge3} and \ref{folge2} now we assume that $G\subset O(p,q)$ is connected and acts irreducibly. If $G$ is not semisimple, its Lie algebra $\lag$ has a non trivial center $\laz$, but the orthogonality of the representation implies that projection of the center on $\rr Id$ is trivial. Hence the  representation is not of real type, i.e. $n=p+q$ is even, and $\laz=\rr J$ where $J$ is the complex structure which commutes with $\lag$. But on the other hand $J\in \laso(p,q)$, i.e. $J$ is compatible with the inner product, which gives that $p$ and $q$ are even as well. Thus, by proposition \ref{Class1}, 
\[\lag\  \subset \ \laso(p,q)\cap \lagl(n/2,\ccc)\ =\ \mf{u}(p/2,q/2).\]
which is the statement of Corollary \ref{folge3}. Furthermore, a straightforward calculation (for conventions see e.g. \cite{baumkath99}) gives that the complex structure $J$ is mapped onto an isomorphism of the spinor module under the spin representation of $\lag$. Hence, if the spinor module of $\lag$ has a trivial submodule, $\lag$ and thus $G$ have to be semisimple. This is the statement of Corollary \ref{folge2}.

\begin{ex}
An example for an irreducible real representation of a Lie group with 2-dimensional center is the
reellification of the representation of $S^1\times CO(n,\rr)$ on $\ccc^n$. The Lie algebra consists of the
matrices $\left(
\begin{array}{cc}
A    & a I_n    \\
-aI_n    &  A
\end{array}
\right)\in \lagl(2n,\rr)$ with $a\in \rr\ \mbox{ and } A\in \mathfrak{co}(n, \mathbb{R})$, where $I_n$ denotes the $n$-dimensional
unit matrix. The center of the identity component of this group is  $S^1\times \rr^+$, the semisimple part is $SO(n)$. In the same manner we can built an example where the center is a spiral in $\ccc^*$  by taking as Lie algebra 
\[ \lag\ :=\ 
\left\{\left.\left(
\begin{array}{cc}
A    & 0     \\
0   &  A
\end{array}
\right)\right| A\in \lason\right\}\ \oplus\  \rr\cdot \left(
\begin{array}{cc}
I_n    & I_n    \\
-I_n   &  I_n
\end{array}\right)\ \subset \ \lagl(2n,\rr),
\]
and as group $G$ the connected subgroup in $Gl(2n,\rr)$ with this Lie algebra. Both groups do not act orthogonally.
\end{ex}

Now we can go ahead with the proof of Theorem \ref{p1}.
Let $G$ be a connected, irreducibly acting Lie subgroup of $Gl(n,\rr)$, and $\lag$
be its Lie algebra. 
Our proof now relies on the following result of \cite{yosida37} and \cite{malcev45} (see also \cite{goto48} where it is a corollary to a deeper result).

\begin{theo}\cite{yosida37}, \cite{malcev45},\cite{goto48} \label{malcev}
\label{f3}
 A connected Lie subgroup  of $Gl(n,\rr)$ is closed in $Gl(n,\rr)$ if and only if its radical is   closed. In particular, if it is semisimple, it is closed.
\end{theo}
Recall that the radical of $G$ is the connected Lie subgroup of $G$ which corresponds to maximal solvable ideal in the Lie algebra $\lag$.
Thus we have to show, that the radical of $G$ is closed in $Gl(n,\rr)$. But by the remarks above, the Lie algebra of $G$ is reductive, and thus
 the radical of $G$ is equal to the identity component of its center, denoted by $Z(G)_0$. Now  the closure $\overline{G}$ of $G$ is still connected, acts irreducibly and has a reductive Lie algebra. By Theorem \ref{malcev} the identity component of its center $Z\left(\ol{G}\right)_0$ is closed in $Gl(n,\rr)$. But  $Z(G)\subset Z(\ol{G})$ because for $z\in Z(G)$ and $g=\lim g_n\in \ol{G}$ it is
\[z\cdot g\ =\  z\cdot \lim g_n\ =\  \lim (z\cdot g_n)\ =\ 0.\]

If we now assume that $G$ is not closed we get by  Theorem \ref{malcev}  that   $Z(G)_0$ is not closed in $Gl(n,\rr)$, i.e.
\[Z(G)_0\ \subsetneqq\  \ol{Z(G)}_0 \subset Z(\ol{G})_0.\]
Now, since $\ol{G}$ is irreducible and connected, Proposition \ref{center} leaves us only with the possibility that  $Z(\ol{G})_0$ is isomorphic  to $\ccc^*$ and $Z(G)_0$ is a one-parameter subgroup of $\ccc^*$. But these are  closed in $\ccc^*$. We obtain
 a contradiction which completes the proof  of Theorem \ref{p1}.

\bigskip

Since holonomy groups are Lie subgroups of $Gl(n,\rr)$, the first point of 
Corollary   \ref{folge1} is a direct consequence of the Theorem \ref{p1}. The second can be obtained by a theorem in \cite{wu67} which contains  several results with different algebraic conditions for subgroups of
the pseudo-orthogonal group, having consequences for holonomy groups.
\btheo \cite{wu67}\label{wutheo}
The following subgroups of $Gl(p+q)$ are closed:
\bnum
\item reductive, indecomposable subgroups of $O(p,q)$,
\item indecomposable subgroups of $O(p,q)$ if $p+q<6$,
\item holonomy groups of affine symmetric spaces.
\enum
\etheo
Here `indecomposable' means `no non-degenerate invariant subspace'. 
One should remark that
the restriction to the dimension in the second point is sharp: In \cite{wu67} is constructed a 6-dimensional K\"ahler manifold
whose reduced holonomy group is non-closed in $SO(4,2)$; also the Lorentzian examples in \cite{bb-ike93} are constructed in dimension 6.
Also in \cite{wu67} is constructed an example of a symmetric
space with solvable, non-Abelian holonomy group which shows that the third point does not follow from the first. Some of these examples are obtained by constructing subgroups containing a torus, which has non-closed 1-parameter subgroups. Our proof shows that such a situation can be excluded if the group acts irreducibly.

In order to obtain the second statement of Corollary \ref{folge1},
note that the first point of Theorem \ref{wutheo} implies that semi-Riemannian holonomy groups which act completely reducibly are closed:  by the de Rham/Wu decomposition theorem \cite{wu64} any 
semi-Riemannian holonomy group is a product of indecomposably acting holonomy groups, but if the group is assumed to act completely reducibly it is reductive and hence closed by the first point of Theorem \ref{wutheo}.
Since the dense line on the Clifford torus provides an example of a  completetly reducibly acting group which is not closed in $Gl(2,\ccc)$, such a result cannot be true for holonomy groups of an arbitrary affine connection due to the result in \cite{hano-ozeki}, that any connected linear Lie group can be obtained as the holonomy group of an affine connection. 

\section{Irreducibly acting, connected subgroups of \mb{$O(1,n)$}}\label{sec3}
In this section we want to give a short proof of Theorem \ref{p2}, 
that the only connected subgroup $G$ of
$O(1,n)$ which acts irreducibly on the Lorentzian space $\R^{1,n}$
is the connected component of the identity of $O(1,n)$ i.e.
$G=SO(1,n)_0$.  This statement was proven in
 \cite{olmos-discala01} where the main goal was to generalize to real
hyperbolic space the following result about minimal homogeneous submanifolds
i.e. orbits of isometry subgroups,  in the Euclidean
space. 
\begin{theo}\cite{D1}
A (extrinsically) homogeneous minimal submanifold of the Euclidean
space must be totally geodesic.
\end{theo}
It turns out that such result also holds in the real hyperbolic
space (see \cite{olmos-discala01} for details).
It is interesting to remark that further investigations of minimal
homogeneous submanifolds were done in several directions
\cite{AD}, \cite{D2}. In particular, the following conjecture was
posed in \cite{D2}.

\begin{con}
Let $M$ be a Riemannian manifold that is either locally
homogeneous or Einstein. Then, any minimal isometric immersion
$f:M \rightarrow \R^n$ must be totally geodesic.
\end{con}

Now, in order to prove  Theorem \ref{p2} we assume that $G\subset O(1,n)$ acts irreducibly and is connected. By Corollary \ref{folge3} it is semisimple and closed by Theorem \ref{p1}.  Our proof requires the following theorem.

\begin{theo}(Karpelevich \cite{Ka}, Mostow \cite{Mo})
Let $M =  {Iso}(M)/K$ be a Riemannian symmetric space of non-compact type.
Then any connected and semisimple subgroup  $G$ of the full isometry
group ${Iso}(M)$ has a totally geodesic orbit $G\cdot p
\subset M$.
\end{theo}

Since $G\subset O(1,n)$ is semisimple this theorem applies to our situation. It implies that the action of $G$ on the hyperbolic
space $H^n = SO(1,n)/SO(n) \subset \R^{1,n}$ is transitive.
Indeed, if the totally geodesic orbit $G\cdot p$ is not the whole
hyperbolic space $H^n$ then $G\cdot p$ is contained in a Lorentzian
subspace $ \mathbb{L}$ of $\R^{1,n}$. This is due to the fact that totally geodesic
submanifolds of $H^n$ are intersections $H^n \cap \mathbb{L}$ where $ \mathbb{L}$
is a Lorentzian subspaces of $\R^{1,n}$. Thus, $G$ can not act
irreducibly as we had assumed. 

Now, 
let $K$ be a maximal connected compact subgroup of the semisimple group $G$. Then by
Cartan's fixed point theorem $K$ has a fixed point $p \in  H^n$.
Since $(G_p)_0$ is compact we get $K = (G_p)_0$. Thus, $(G,K)$ is
a symmetric pair such that $H^N = G/K$. Then, from the uniqueness of
such symmetric pairs
(see \cite[pp. 243]{Helgason78}) we get $G= SO(1,n)_0$ and $K=SO(n)$. This proves Theorem \ref{p2}.

\bigskip
 
A different,  almost algebra-free
proof of Theorem \ref{p2} which uses dynamical methods, can be found in \cite{BZ}.
A purely algebraic proof was given in \cite{benoist-harpe03}.

\section{Invariant bilinear forms of irreducible representations}\label{sec4}



As in the second section we consider an irreducibly acting subgroup $G \subset Gl(V)$ of a real vector space $V$ and denote by
$B_G(V)$ the vector space of $G$-invariant bilinear forms. If $B_G(V)$ is non-trivial it is intimately connected to $End_G(V)$. By  Schur's lemma
 a non-zero $a \in B_G(V)$ is non-degenerate since its kernel is $G$-invariant and not equal to $V$.
Thus Riesz' theorem provides a one-to-one correspondence between $B_G(V)$ and $End_G(V)$ via
\beq\label{Riesz}
\rca{1.5}{
\ba{rcl} End_G(V) &\overset{\sim}{\rightarrow}&  B_G(V) \\
B &\mapsto& b=a(B(\cdot),\cdot). 
\ea}
\eeq
In particular this map endows $B_G(V)$ with the structure of an associative algebra and by Proposition \ref{End_G(V)}  $B_G(V)$ is isomorphic to $\R$, $\C$, or $\HH$.
The unique decomposition of a bilinear form into symmetric and skew-symmetric parts applies
also to $G$-invariant bilinear forms, since the (skew-)symmetrization of a $G$-invariant form inherits this property. Thus
\beq\label{decomp1}
B_G(V)=S_G(V) \oplus \Lambda_G(V).
\eeq
This induces a decomposition of $End_G(V)$ into $a$-selfadjoint and $a$-skewadjoint operators,
\beq\label{decomp2}
End_G(V)=S_G^a(V) \oplus \Lambda_G^a(V).
\eeq
If $a$ is symmetric, $S_G(V)$ corresponds to $S_G^a(V)$ under (\ref{Riesz}) and
if $a$ is skew-symmetric, $S_G(V)$ corresponds to $\Lambda_G^a(V)$.
The main question is what are the possible dimensions of $S_G(V)$ and what are
the occurring signatures. A first answer gives the following statement.

\begin{satz}\label{adjointness}
Let $a,b \in B_G(V)$ be linear independent,
$b=a(B(\cdot,\cdot))$ by {\rm (\ref{Riesz})} and
$B=\alpha \,Id + \beta J$ according to Corollary  {\rm \ref{id+J}}.
\bnum
\item[(i)] If $a$ and $b$ are both symmetric (or skew-symmetric),
$J \in S^a_G(V)$ and thus $J$ is an anti-isometry with respect to both $a$ and $b$.
In particular $\Sig(a)=\Sig(b)=(n/2,n/2)$ where $n=\dim V$.

\item[(ii)] If $a$ is symmetric and $b$ is skew-symmetric,
$B = \beta J \in \Lambda^a_G(V)$
and thus $J$ is an isometry with respect to both $a$ and $b$.
\enum
\end{satz}

\begin{proof}
(i) Since $B$ and $Id$ are $a$-selfadjoint, the same holds for $J$.
Using $[B,J]=0$ we obtain
\begin{gather*}
a(J(x),J(y))=a(J^2(x),y)=-a(x,y) \text{, and}\\
b(J(x),J(y))=a(B \circ J(x),J(y))=-a(B(x),y)=-b(x,y).
\end{gather*}

(ii) Here $B$ is $a$-skewadjoint. This implies for its minimal polynomial
$\mu_B(x)=\mu_{-B}(x)=\mu_B(-x)$, hence $\mu_B(x)=x^2+\beta^2$ (cf. Corollary  \ref{id+J}),
i.e. $B=\beta J$. The remaining part is analogous to (i).
\end{proof}

Next we determine all possible pairs
$\big(\dim S_G(V),\dim \Lambda_G(V)\big)$ by describing their maximal
representations analogous to Proposition \ref{Class1} of Section \ref{sec1}. Recall that a representation is self-dual if and only if the space of non-degenerate invariant bilinear forms is non-trivial.
\begin{satz}\label{Class2}
Let $\kappa:G \rightarrow Gl(n,\R)$ be an irreducible self-dual representation on $\R^n$.
Then up to conjugation $\kappa(G)$ is contained in one of the following subgroups $L \subset Gl(n,\R)$ with $p+q=n$:
\[
\rca{2}{
\ba{c|c|c|l}
End_{\kappa(G)}(\R^n) & \dim S_{\kappa(G)}(\R^n)& \dim \Lambda_{\kappa(G)}(\R^n)& L \subset Gl(n,\R) \\ \hline
\R & 1 & 0 & O(p,q) \\ 
   & 0 & 1 & Sp(n/2,\R) \\ \hline
\C & 2 & 0 & O(n/2,\C) ,\;{\sst(n \geq 4)}\\ 
   & 0 & 2 & Sp(n/4,\C) \\  
   & 1 & 1 & U(p/2,q/2)  \\ \hline
\HH & 1 & 3 & Sp(p/4,q/4),\;{\sst(n \geq 8)} \\ 
    & 3 & 1 & O^*(n/4),\;{\sst(n \geq 8)} \\
\ea
}
\]
Moreover we have the following isomorphisms\footnote{
For details how the groups are embedded into $Gl(n,\R)$ resp. $Gl(n,\C)$ we refer to the proof.}
\beq\label{identifications}
\rca{1.5}{\ba{rcl}
O(n/2,\C) &\simeq& O(n/2,n/2) \cap Gl(n/2,\C) \\
Sp(n/2,\C) &\simeq& Sp(n,\R) \cap Gl(n/2,\C) \\
U(p/2,q/2) &\simeq& O(p,q) \cap Sp(n,\R) \\
Sp(p/4,q/4) &\simeq& U(p/2,q/2) \cap Sp(n/4,\C)  \\
O^*(n/4) &\simeq& U(n/4,n/4) \cap O(n/2,\C).
\ea}
\eeq
\end{satz}

\begin{bem}
This proposition raises the question if there are proper subrepresentations of the different groups $L \subset Gl(n,\R)$.
The answer depends very much on the group $L$ in question. To illustrate this, note that any compact simple Lie group admits
irreducible representations in arbitrary high dimensions. All these representations are contained in $O(n)$
due to Weyl's trick. Considered the other way around $O(n)$ has in general a lot of irreducible subrepresentations.
In contrast, there are no proper subgroups of $SO_0(1,n)$ which act irreducibly, see section \ref{sec3}.

\end{bem}

\begin{proof}
For general considerations set $V:=\R^n$; we will return at the end to $\R^n$ by choosing an appropriate basis.
First note that $(\kappa,V)$ is self-dual if and only if $B_{\kappa(G)}(V)\neq \{0\}$.
In particular $End_{\kappa(G)}(V) \simeq B_{\kappa(G)}(V)$ according to (\ref{Riesz})
and we may distinguish between the various types of $(\kappa,V)$. We determine in each case the maximal subgroup $L \subset Gl(V)$
fixing every element of $B_{\kappa(G)}(V)$ and thus $\kappa(G)\subset L$.

\vspace{1em}
\mb{$(\kappa,V)$ of real type:}
This is the simplest case, since $B_{\kappa(G)}(V)$ is 1-dimensional
and thus spanned either by a symmetric or skew-symmetric bilinear form $a$ (cf. (\ref{decomp1})).
In the symmetric case we can find a (pseudo-)orthonormal basis, i.e.
\[
(a_{ij}) = \II_{p,q} := \lp \ba{cc} -\II_p & \\ & \II_q \ea \rp.
\]
Its isometry group is the (pseudo-)orthogonal group
\[
O(p,q):=\menge{A \in Gl(n,\R)}{A^t \cdot \II_{p,q} \cdot A = \II_{p,q} }.
\]
In the skew-symmetric case we can find a symplectic basis, i.e.
\[
(a_{ij})= \JJ_{n/2}:= \lp \ba{cc} & -\II_{n/2} \\ \II_{n/2} \ea\rp.
\]
Its isometry group is the real symplectic group
\[
Sp(n/2,\R):=\menge{ A \in Gl(n,\R) }{ A^t \cdot \JJ_{n/2} \cdot A = \JJ_{n/2} }.
\]

\vspace{1em}
\mb{$(\kappa,V)$ of complex type:}
First we fix a (skew-)symmetric $a \in B_G(V)$ and consider its bilinear extension $a_{\C}$
as well as its sesquilinear extension $\ol{a}_{\C}$ to the complexification $V_{\C}$:
\[
\rca{1.5}{
\ba{lcr}
a_{\C}(x+iy,u+iz)&:=&a(x,u)-a(y,z) + i \big( a(x,z) + a(y,u) \big), \\
\ol{a}_{\C}(x+iy,u+iz)&:=&a(x,u)+a(y,z) + i \big( a(x,z) - a(y,u) \big). \ea
}
\]
$a_{\C}$ is (skew-)symmetric, $\ol{a}_{\C}$ is (skew-)Hermitian and they are linked by
the formula $a_{\C}(\ol{v},w)=\ol{a}_{\C}(v,w)$. From this follows
\beq\label{a_C}
 a_{\C}(\ol{v},\ol{w}) = \ol{a_{\C}(v,w)} \qquad \text{ and } \qquad
\ol{a}_{\C}(\ol{v},\ol{w}) = \ol{\ol{a}_{\C}(v,w)}.
\eeq
Exactly one of them has to vanish on the $\kappa_{\C}$-irreducible subspace $W$.
Indeed, if we suppose both to be non-zero and set $a_{\C}\big|_W=\ol{a}_{\C}\big|_W(J(\cdot),\cdot)$
one shows that $J \in \ol{End_{\rho}(W)}$, hence $J=0$.
On the other hand $a_{\C}\big|_{W \times W} = 0$ together with
(\ref{a_C}) implies $a_{\C}\big|_{\ol{W} \times \ol{W}} = 0$.
Thus $a_{\C}\big|_{\ol{W} \times W} = \ol{a}_{\C}\big|_{W \times W}$ has to be non-degenerate
and vice versa. Lets denote by $\wt{a}$ the non-vanishing form on $W$.
Since $\wt{a}$ is $\rho$-invariant it induces the $\kappa$-invariant
$\C$-valued $\R$-bilinear form $\psi_* \wt{a}$ on $V$ via (\ref{psi}).
In the following we will suppress the isomorphism $\psi$.
Real and imaginary part of $\wt{a}$ are related by
\[
Im \, (\wt{a})(x,y) = - Re\, (\wt{a})(x,I(y)).
\]
In particular they are linear independent and thus
$B_{\kappa(G)}(V)$ is spanned by these two forms. So the isometry group
of $\wt{a}$ is isomorphic to the maximal subgroup of $L \subset Gl(V)$ which fixes
any element of $B_{\kappa(G)}(V)$. Note that any element of $L$ commutes with $I$ and thus it is complex linear.

If $\wt{a}$ is symmetric, the same holds for its real and imaginary part and their
signature has to be $(n/2,n/2)$ (cf. Proposition \ref{adjointness}(i)).
We can find a complex orthonormal basis, i.e. $(\wt{a}_{ij})=\II_{n/2}$ and
the isometry group is the complex orthogonal group
\[
O(n/2,\C)=\menge{A \in Gl(n/2,\C)}{A^t \cdot A = \II_{n/2}}.
\]

If $\wt{a}$ is skew-symmetric, the same holds for its real and imaginary part.
We can find a complex symplectic basis i.e. $(\wt{a}_{ij})=\JJ_{n/4}$ and
the isometry group is the complex symplectic group
\[
Sp(n/2,\C)=\menge{A \in Gl(n/2,\C)}{A^t \cdot \JJ_{n/4} \cdot A = \JJ_{n/4}}.
\]

Finally $\wt{a}$ might be Hermitian (a complex skew-Hermitian form turns into a Hermitian one
by multiplication with $i$). We can find a complex (pseudo-)orthonormal basis, i.e.
$(\wt{a}_{ij})=\II_{p/2,q/2}$ and the isometry group is the unitary group
\[
U(p/2,q/2):=\menge{A \in Gl(n/2,\C)}{\ol{A}^t \cdot \II_{p/2,q/2} \cdot A = \II_{p/2,q/2} }.
\]
In this case the real part is symmetric and has signature $(p,q)$ and the imaginary
part is skew-symmetric.

As mentioned in Proposition \ref{Class1} the complex basis $\{b_i\}_{i=1}^{n/2}$ of $V$ (with respect to $I$)
induces the real basis $\{b_i, I(b_i)\}_{i=1}^{n/2}$.
Thus $I=\JJ_{n/2}$ and we obtain the embedding
\[
\rca{1.5}{
\ba{rcl}
Gl(n/2,\C) &\simeq& \menge{ C \in Gl(n,\R) }{ C \circ \JJ_{n/2}=\JJ_{n/2} \circ C }\\
A+iB  &\mapsto& \lp \rca{1}{\ba{cc} A& -B \\ B& A \ea} \rp.
\ea}
\]
Now the real part of the symmetric form $\wt{a}=\II_{n/2}$ is given in the associated real basis
by $Re (\wt{a})=\lp \bsm \II_{n/2}& \\ &  -\II_{n/2}\esm\rp $.
Its isometry group is $O(n/2,n/2)$, thus we get the first identity of (\ref{identifications}).
Analogous, the real part of the symplectic form $\wt{a}=\JJ_{n/4}$ is given by
$Re (\wt{a})= \lp \bsm \JJ_{n/4}& \\ &  -\JJ_{n/4}\esm\rp $.
Its isometry group is conjugated to $O(n/2,n/2)$ which yields the second identity of (\ref{identifications}).
The real part of the Hermitian form $\wt{a}=\II_{p/2,q/2}$ is given by
$Re (\wt{a})= \lp \bsm \II_{p/2,q/2}& \\ & \II_{p/2,q/2} \esm\rp $. Its isometry group is conjugate to $O(p,q)$.
Instead of taking the intersection with the centralizer of $\JJ_{n/2}$ as above we take here the isometry group of
the imaginary part $Im (\wt{a})= \JJ_{n/2} $ which is $Sp(n/2,\R)$, hence the third identity of (\ref{identifications}).

\vspace{1em}
\mb{$(\kappa,V)$ of quaternionic type:}
For representations of real or complex type all possible dimensions for the
subspaces $S_{\kappa(G)}(V)$ and $\Lambda_{\kappa(G)}(V)$ occurred.
This is no longer true in the quaternionic case.

\begin{lem}\label{Lemma}
If $(\kappa,V)$ is self-dual and of quaternionic type then
$S_{\kappa(G)}(V)$ and  $ \Lambda_{\kappa(G)}(V)$ are odd-dimensional, i.e. their dimension
is 1 and 3. In particular, $\kappa$ is both, orthogonal and symplectic.
If the 1-dimensional subspace is spanned by $\{a\}$ then under the identification
$End_{\kappa(G)}(V) \simeq \HH$ the decomposition {\rm (\ref{decomp2})} is given by
\[
Re (\HH) = S^a_{\kappa(G)}(V) \qquad Im(\HH)=\Lambda^a_{\kappa(G)}(V).
\]
\end{lem}

\begin{proof}
Clearly $Re(\HH)=\R \cdot Id \subset S^a_{\kappa(G)}(V)$. On the other hand
$\Lambda^a_{\kappa(G)}(V) \subset Im(\HH)$ by Proposition \ref{adjointness}, hence
$Im(\HH)= \big( S^a_{\kappa(G)}(V) \cap Im(\HH) \big) \oplus \Lambda^a_{\kappa(G)}(V)$.
One of the subspaces has dimension greater or equal than two and is spanned by
anti-commuting complex structures $I,J$. Irrespective of whether $I,J$ are
self- or skewadjoint with respect to $a$, their product is skewadjoint:
$(I \circ J)^*=J^* \circ I^* = J \circ I = - I \circ J$.
\end{proof}

As in Proposition \ref{Class1} we consider $V$ as right $\HH$-vector space
via $x\cdot \lambda=\ol{\lambda}(x)$.
Then an element $a \in B_{\kappa(G)}(V)$ as in the preceeding lemma yields
the following quaternionic sesquilinear form on $V$:
\[
a_{\HH}(x,y):=a(x,y)+i\cdot a(xi,y)+j\cdot a(xj,y)+k\cdot a(xk,y).
\]
Recall that one has to check $a_{\HH}(x\lambda,y)=\ol{\lambda}a_{\HH}(x,y)$
and $a_{\HH}(x,y\lambda)=a_{\HH}(x,y)\lambda$. Since multiplication (from the right)
with imaginary quaternions is an $a$-skew\-adjoint operation according to Lemma \ref{Lemma}, $a_{\HH}$ is
Hermitian if $a$ is symmetric and skew-Hermitian otherwise.

By construction, $B_{\kappa(G)}(V)$ is spanned by the real and imaginary parts of $a_{\HH}$,
hence the group $L$ which fixes all elements of $B_{\kappa(G)}(V)$ is the isometry group of $a_{\HH}$
(again $\HH$-linearity is ensured already by leaving $a_{\HH}$ invariant).

Now, for any (skew-)Hermitian
form one can find an orthogonal basis \cite[Ch. I, \S 8]{Dieudonne71}.
In the Hermitian case the basis can be normed to the length $\pm 1$. Thus the isometry
group is the quaternionic unitary group
\[
Sp(p/4,q/4)=\menge{A\in Gl(n/4,\HH)}{\ol{A}^t \cdot \II_{p/4,q/4} \cdot A = \II_{p/4,q/4} }.
\]
In the skew-Hermitian case the basis can be normed to the length $i$. Thus the isometry
group is
\[
O^*(n/4)=\menge{A \in Gl(n/4,\HH)}{\ol{A}^t \cdot i\II_{n/4} \cdot A = i \II_{n/4}}.
\]
The embedding $L \subset Gl(n,\R)$ follows from the embedding $Gl(n/4,\HH) \subset Gl(n,\R)$ as in Proposition \ref{Class1}.
In order to obtain the remaining identities of (\ref{identifications}), we fix $i$ as complex structure and thus
represent any quaternionic matrix by two complex matrices: $A+iB+jC+kD=(A+iB)+(C+iD)j=U+Vj$.
This yields the algebra isomorphism
\[
\rca{1.5}{\ba{rcl}
Gl(n/4,\HH) &\simeq& \menge{C \in Gl(n/2,\C)}{ \ol{C} \circ \JJ_{n/4} = \JJ_{n/4} \circ C} \\
U+Vj &\mapsto& \lp \rca{1}{\ba{cc} U &-V \\ \ol{V} & \ol{U} \ea}\rp.
\ea}
\]
Since under this identification the operation $C \mapsto \ol{C}^t$ is the same in $Gl(n/4,\HH)$ and $Gl(n/2,\C)$ it is easily seen, that
$Sp(p/4,q/4)$ is equal to the intersection of the isometry group $U(p/2,q/2)$ of the Hermitian form
$\lp \bsm \II_{p/4,q/4}&\\&\II_{p/4,q/4}\esm\rp$
with the isometry group $Sp(n/4,\C)$ of the symplectic form  $\lp \bsm &-\II_{p/4,q/4}\\\II_{p/4,q/4}&\esm\rp$.
Analogously we obtain $O^*(n/4)$ as intersection of the isometry group $U(n/4,n/4)$ of the skew-Hermitian form
$\lp \bsm i\II_{n/4}& \\ &-i\II_{n/4} \esm\rp$ with the isometry group $O(n/2,\C)$ of the symmetric form
$\lp \bsm &\II_{n/4}\\\II_{n/4}\esm\rp$. This yields the remaining identities of (\ref{identifications}).

We conclude the proof by showing that the maximal groups $L \subset Gl(n,\R)$ are acting irreducibly on $\R^n$.
One knows even more: For any subgroup $L \subset Gl(n/d,\K) \subset Gl(n,\R)$ occuring in the list of the proposition
its centralizer coincides with the corresponding homothety group $H_{\K}$:
\[
End_L(\R^n)=H_{\K},\quad L \subset Gl(\K).
\]
For the symplectic groups this can be easily verified.
For the unitary groups this is true beginning with $n/d \geq 2$ and for the orthogonal groups
it is true for $n/d \geq 3$ (see \cite[Ch. II, \S 3]{Dieudonne71}).
In this context the quaternionic groups $Sp(p/4,q/4)$ and $O^*(n/4)$ are comprehended as unitary groups.
Since any homothety is invertible the above groups act irreducibly, otherwise the projection onto an invariant subspace would be an
element of the centralizer which is certainly not invertible. It remains to discuss irreducibility in the excluded small dimensions.
\end{proof}

\begin{bem}
A quaternionic vector space does not admit any symmetric or skew-symmetric bilinear form. This is reflected
in the fact that the space of symmetric or skew-symmetric bilinear forms is never 4-dimensional (cf. Lemma \ref{Lemma}).
\end{bem}

\begin{bem}
Changing the complex basis by the matrix $\frac{1}{\sqrt{2}}\lp\bsm 1 & i \\ 1 & -i \esm\rp$
we obtain the embedding $O^*(n/4) \subset Gl(n/2,\C)$ as given in \cite[Ch.X,\S 2,1.]{Helgason01}.
There it is explained how $O^*(n/4)$ occurs as dual of the symmetric space $O(n/2) \big/ U(n/4)$ which justifies the notation.
\end{bem}

The considerations above can be generalized to non-irreducible representations of Lie groups or Lie algebras.
This has been done in \cite{MedinaR93}. Of course the structure of the algebra $End_G(V)$ becomes more involved.
On the other hand we may restrict our attention to special representations as e.g. the adjoint representation $Ad(G)\subset Gl(\g)$
of a Lie group $G$.
To ask for an $Ad(G)$-invariant non-degenerate symmetric bilinear form on $\g$ becomes interesting from a
geometrical point of view, since any such bilinear form induces a pseudo-Riemannian metric on $G$ which is invariant under left-
and right-multiplication. In particular $G$ becomes a symmetric space. Hereafter we cite some results in this direction.

As shown above there are representations which are symplectic but not orthogonal. This fails for adjoint representations:
\begin{satz}[\cite{MedinaR93}, Theorem 1.4]
A Lie algebra $\g$ admits an $\ad(\g)$-invariant non-degenerate symmetric bilinear form if and only if it is self-dual.
\end{satz}

On the other hand it has been shown:
\begin{satz}[\cite{MedinaR93}, Corollary  1.7]
A Lie algebra $\g$ admits an $\ad(\g)$-invariant skew-symmetric bilinear form if and only if $\codim_{\g}[\g,\g]\geq 2$.
\end{satz}

In particular, for simple Lie algebras the adjoint representation is irreducible and $[\g,\g]=\g$. Thus,
 they cannot be symplectic which excludes many cases of Proposition \ref{Class2}.


\end{document}